\documentclass[12pt]{iopart}

\usepackage{amsmath}
\usepackage{amsfonts}
\usepackage{amsthm}
\usepackage{amssymb}
\usepackage{color}
\usepackage{pictex}
\usepackage{enumerate}
\usepackage{cite}

\newtheoremstyle{thm}
  {\baselineskip}
  {\baselineskip}
  {\itshape}
  {}
  {\bf}
  {.}
  {.5em}
  {}
  
\newtheoremstyle{others}
  {\baselineskip}
  {\baselineskip}
  {\upshape}
  {}
  {\bf}
  {.}
  {.5em}
  {}

\theoremstyle{thm}
\newtheorem{theorem}{Theorem}[section]
\newtheorem{proposition}[theorem]{Proposition}

\newtheorem{corollar}[theorem]{Corollary}

\theoremstyle{others}
\newtheorem{counterexample}[theorem]{Counterexample}
\newtheorem{example}[theorem]{Example}

\newtheorem{remark}[theorem]{Remark}

\newcommand{\argmin}{\operatorname{argmin}}
\newcommand{\besov}{\mathcal{B}}

\newcommand{\ip}[2]{\langle #1, #2 \rangle}
\newcommand{\norm}[1]{\|#1\|}

\newcommand{\R}{\mathbb{R}}
\newcommand{\rg}{\operatorname{rg}}
\newcommand{\set}[2]{\{#1\, | \,#2\}}

\newcommand{\sign}{\operatorname{sign}}
\newcommand{\bigO}{\mathcal{O}}

\newcommand{\Z}{\mathbb{Z}}

\begin{document}

\title[Regularization in Besov scales]{Optimal Convergence Rates for Tikhonov Regularization in Besov Scales}

\author{D A Lorenz$^1$ and D Trede$^1$}

\address{$^1$Zentrum f\"ur Technomathematik, University of Bremen, D--28334 Bremen, Germany}

\eads{\mailto{dlorenz@math.uni-bremen.de}, \mailto{trede@math.uni-bremen.de}}

\begin{abstract}
  In this paper we deal with linear inverse problems
  and convergence rates for Tikhonov regularization.
  We consider regularization in a scale of Banach spaces,
  namely the scale of Besov spaces.
  We show that regularization in Banach scales differs
  from regularization in Hilbert scales in the sense that it is possible that stronger source conditions may lead to weaker convergence rates and vive versa.
  Moreover, we present optimal source conditions for regularization in Besov scales.
\end{abstract}

\section{Introduction}
Regularization of inverse problems formulated in Banach spaces have been of recent interest.
On the one hand there are a several theoretical regularization results such as convergence rates in a general Banach spaces setting,
see e.g.~\cite{burger2004convarreg,resmerita2005regbanspaces,resmerita2006nonquadreg,
burger2007errorinvscalespace,hofmann2007convtikban,
hein2008convratesbanach},
or convergence rates for special sequence spaces such as $\ell^p$, $1 \leq p <2$, 
see e.g.~\cite{lorenz2008reglp,ramlau2008regproptikhonov}.
On the other hand there are results which deal with solving inverse problems formulated in Banach spaces,
such as Landweber-like iterations or minimization methods for Tikhonov functionals,
see e.g.~\cite{schoepfer2006illposedbanach} and~\cite{daubechies2003iteratethresh,bredies2008harditer,bredies2008itersoftconvlinear,bredies2008forwardbackwardbanach,bonesky2007gencondgradnonlin,ramlau2006tikhproject,griesse2008ssnsparsity,schuster2008tikbanach},
respectively.

The interest in Banach spaces is due to the fact that in many situations a Banach space is better suited to model the data under consideration than a Hilbert space.
In the context of image processing, for example, the Banach space $BV$ of functions of bounded variation is used to model images with discontinuities along lines~\cite{rudin1992tv,acar1994analysis,burger2004convarreg}.
Moreover,~\cite{hofmann2007convtikban} presents two examples in which the use of Banach spaces is necessary for a thorough formulation of the problem.
Another class of Banach spaces are the Besov spaces $B^s_{p,q}$ which play an important role in inverse problems related to image processing, see e.g.~\cite{chambolle1998shrink,chan2006imageprocessing,lorenz2007wavproj}.

In this paper we make a first attempt to analyze inverse problems in scales of Banach spaces generalizing classical Hilbert scales~\cite{natterer1984hilbertscales}.
The easiest scale of Banach spaces is the scale of Sobolev spaces $W^s_p$.
However, we are going to use Besov spaces $B^s_{p,p}$ since they coincide with the Sobolev scale in most cases if the integrability indices coincide.
Moreover, they come with a characterization in terms of wavelet coefficients which make them easy to use for our purposes.
We apply previous convergence rates regularization results from~\cite{burger2004convarreg,hofmann2007convtikban} in the scale of Besov spaces and develop optimal convergence rates.
To this end, we derive the source conditions that lead to a convergence rate of $\bigO(\sqrt{\delta})$ in a certain Sobolev space.

Consider the equation
\begin{equation}\label{op_eq}
  F \; u^\dagger = v,
\end{equation}
where $F$ is a linear continuous operator
\[
  F: \besov_D \to L_2
\]
between the Besov space $\besov_D:=B_{p_D,p_D}^{s_D}$,
$s_D \in \R$, $p_D > 1$,
and Lebesgue space $L_2$.
In general these function spaces contain functions or distributions defined on the subset $\Omega \subset \R^d$.
Due to clarity we omit $\Omega$ in the following.
The Besov spaces $B_{p,p}^s$ are subspaces of the
space of tempered distributions $\mathcal{S}'$ and, in contrast to $\mathcal{S}'$,
they are Banach spaces for $p\geq 1$~\cite{runst1996spaces}.
Different from classical approaches we use the domain $\besov_D$%
---often a superset of $L_2$---and not $L_2$ itself.
That may be of interest in some applications, e.g. mass-spectrometry where the data consists of delta peaks
(see \cite{klann2007shrink_vs_deconv,dahlke200multiscale}) which are not elements of $L_2$.

If we assume that only noisy data $v^\delta$ with noise level $\norm{v-v^\delta} \leq \delta$ are available,
the solution of (\ref{op_eq})
could be unstable and has to be stabilized by regularization methods.
We use regularization with a Besov constraint, i.e.~we regularize by minimizing a not necessarily quadratic functional $T_\alpha:\besov_D\to [0,\infty]$ defined by
\begin{equation}\label{tikhonov_eq}
  T_\alpha(u) := \norm{Fu - v^\delta}_{L_2}^2 + \alpha \norm{u}_{\besov_R}^{p_R},
\end{equation}
where $\besov_R:=B_{p_R,p_R}^{s_R}$ is a Besov space,
not necessarily equal to $\besov_D$.
Since $T_\alpha$ shall be defined on $\besov_D$ we define
\[
  \norm{u}_{\besov_R} = \infty
\]
for $u\notin \besov_R$.

In this paper we will investigate regularization properties and convergence rates
of the regularization method consisting of the minimization of (\ref{tikhonov_eq}),
i.e. $u^{\alpha,\delta}\in\argmin T_\alpha(u)$.
The proceeding is as follows.
\begin{enumerate}[i)]
  \item In section~\ref{sec_notation} we introduce the notation and collect preliminary results. 
  \item In section~\ref{sec_convergence} we apply convergence rates results for Banach 
    spaces~~\cite{burger2004convarreg,hofmann2007convtikban}.
    With the constraints on $p_R$ and $s_R$ in mind and the
    parameter rule $\alpha = \delta$ we will
    get a stable approximation, i.e.
    \[
      \norm{u^{\alpha,\delta} - u^\dagger}_{\besov_R} \to 0, \quad \delta \to 0,
    \]
    and a convergence rate in the Sobolev space $H^{\sigma}$
    \[
      \norm{u^{\alpha,\delta} - u^\dagger}_{H^{\sigma}}
        =
        \bigO(\sqrt{\delta}),
    \]
    with $\sigma$ depending on $s_R$ and $p_R$ 
    (theorems~\ref{convergence} and~\ref{convergence_rate}).
    These results restrict the choice of possible regularization spaces $\besov_R$.
    Using Besov space embeddings in section~\ref{sec_optimization} we will get a generalization of the first result.
    We find a convergence rate---also formulated in a Sobolev space---which 
    holds for a larger set of Besov space penalties $\norm{\cdot}_{\besov_R}^{p_R}$
    (theorem~\ref{convergence_rate_weaker_source}).
  \item The convergence result gets stronger as $\sigma$ increases since for $\theta >0$ it holds $H^{\sigma+\theta} \subset H^{\sigma}$.  
    Since $\sigma$ depends on $s_R$ and $p_R$, we address the question
    how to choose $\besov_R$ in a way such that $\sigma$ is maximal.
    We will find the regularization penalty $\norm{\cdot}_{\besov_R}^{p_R}$,
    which gives the best estimate with respect to $\sigma$.
  \item In section~\ref{sec_examples} we apply these results to some operators defined in
    Sobolev and Besov spaces to demonstrate the differences.
\end{enumerate}

\section{Notation and Basic Besov Space Properties}
 \label{sec_notation}
As already mentioned, Besov spaces $B_{p,q}^s$ are subspaces of the space of tempered distributions $\mathcal{S}'$.
They coincide with special cases of traditional smoothness function spaces
such as H\"older and Sobolev spaces.
Note e.g., that $B^s_{p,p} = W^s_p$ for $s\notin\Z$ and $p \geq 1$ and even
$B^s_{2,2} = H^s := W^s_2$ for all $s\in\R$.
As from now we use the term Sobolev space only for the Hilbert spaces $H^s$.
This clarifies the characterization that the Besov space $B_{p,q}^s$ contains functions
having $s$ derivatives in $L_p$ norm. The second integrability index $q$ declares a finer nuance of smoothness.
In the following we omit the second integrability index $q$ of the Besov spaces
which is always equal to the corresponding first one $p$.

There are a several ways of defining Besov spaces. Most commonly they are defined via the modulus of smoothness,
a way to model differential properties.
For a detailed introduction of Besov spaces via moduli of smoothness 
in conjunction with other smoothness spaces see e.g. \cite[section 4.5]{devore1998nonapprox}.

Another way defining Besov spaces is based on wavelet coefficients.
According to \cite{cohen2003numwav}
for all $s\in\R$,$\ p>0$,
there exists a wavelet basis 
$\{\psi_\lambda\}_{\lambda\in\Lambda}$
such that
\begin{equation}
  \label{eq_equiv_besov_norm}
\norm{u}_{B^s_p}^{p} \asymp
\sum\limits_{\lambda\in\Lambda} 
2^{p(s+d(\frac{1}{2}-\frac{1}{p})) \; |\lambda|}
|u_\lambda|^{p},
\end{equation}
where $u_\lambda = \ip{u}{\psi_\lambda} = \int u\, \psi_\lambda dx$ are the wavelet coefficients of $u$.
The notation $A\asymp B$ means that there exist constants $c,C>0$ such that $cA\leq B\leq CA$.
We will use this equivalent norm throughout the paper.

An important ingredient in the analysis of the regularization method (\ref{tikhonov_eq}) is the embedding result (cf.~\cite{runst1996spaces}):

\begin{proposition}\label{embed_besov_spaces}
  Let $B_{p_1}^{s_1}$, $B_{p_2}^{s_2}$ be Besov spaces. If
  \begin{equation}\label{eq-embed_besov_spaces}
    s_1 - \frac{d}{p_1} > s_2 - \frac{d}{p_2}
    \quad \mbox{and} \quad
    p_1 \leq p_2,
  \end{equation}
  then $B_{p_1}^{s_1} \subset B_{p_2}^{s_2}$ continuously.
  The term $s - \frac{d}{p}$ is called \emph{differential dimension} of $B_p^s$.
\end{proposition}

\begin{figure}[h]
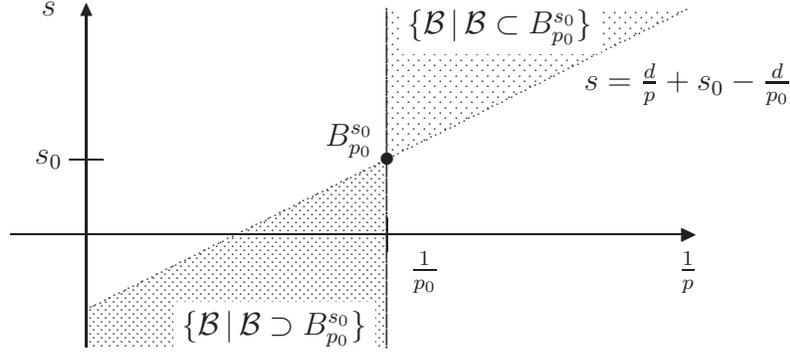

  \[
  \beginpicture
    \setcoordinatesystem units <2cm,1cm>
    \setplotarea x from -0.5 to 4, y from -1.5 to 3
    \inboundscheckon
    \axis left shiftedto x=0
          ticks long in unlabeled at 1 /
                long out unlabeled at 1 /
    /
    \axis bottom shiftedto y=0
          ticks long in unlabeled at 2 /
                long out unlabeled at 2 /
    /
    \put{$\scriptstyle \blacktriangle$} at 0 3
    \put{$s$} at -0.25 3
    \put{$\scriptstyle \blacktriangleright$} at 4 0
    \put{$\frac{1}{p}$} at 4 -0.5
    
    \put{$\bullet$} at 2 1
    \put{$B_{p_0}^{s_0}$} at  1.75 1.25
    
    \put{$s_0$} at -0.25 1
    \put{$\frac{1}{p_0}$} at 2.25 -0.5
    
    \setdots <2pt>
    \plot 0 -1  4 3 /
    \put{$s=\frac{d}{p} + s_0 - \frac{d}{p_0}$} at 4 2
    \setsolid
    \plot 2 -1.5  2 3 /
    
    \setshadegrid span <1.5pt>
    \vshade 0 -1.5 -1
            <z,z,z,z>
            2 -1.5  1
    /
    \put{\colorbox{white}{$\set{\besov}{\besov \supset B_{p_0}^{s_0}}$}} at 1.25 -1.25
    \setshadegrid span <2pt>
    \vshade 2 1 3
            <z,z,z,z>
            4 3 3
    /
    \put{\colorbox{white}{$\set{\besov}{\besov \subset B_{p_0}^{s_0}}$}} at 2.75 2.75
  \endpicture
  \]
  \caption{DeVore diagram for the embedding of Besov spaces.}
  \label{bild_embed_besov_spaces}
\end{figure}

The embedding of Besov spaces is often visualized with the help of the DeVore diagram~\cite{devore1998nonapprox} where one plots the smoothness $s$ against $1/p$, see figure~\ref{bild_embed_besov_spaces}.
By $B_{p_1}^{s_1} \subset B_{p_2}^{s_2}$, in the following,
we denote not only the set-theoretical embedding but also the continuous embedding.

We are going to use the following Besov spaces:
\begin{itemize}
\item $\besov_D := B^{s_D}_{p_D}$ for the domain of $F$.
\item $\besov_R := B^{s_R}_{p_R}$ for the space in which we regularize.
\item $\besov_S := B^{s_S}_{p_S}$ for the source condition.
\item $\besov_G := B^{s_G}_{p_G}$ for the range of $F^*$ which models the smoothing properties of $F$.
\end{itemize}

As stated above, the smoothing properties of the operator
\[
  F: \besov_D \to L_2
\]
are modeled by assuming that the range of its adjoint is small, namely
\[
\rg F^* = \besov_G \subset \besov_D^\ast = B_{p_D^\ast}^{-s_D},
\]
where $p^\ast$ is defined via $\frac{1}{p} + \frac{1}{p^\ast} = 1$,
hence $p^\ast = \frac{p}{p-1}$.
Consequently, we have
\[
  s_G - \frac{d}{p_G} > -s_D - \frac{d}{p_D^\ast},
  \quad
  \frac{1}{p_G} \geq \frac{1}{p_D^\ast}.
\]


\section{Convergence and Convergence Rates}
\label{sec_convergence}
The first result we need is a regularization result in the regularization space $\besov_R$. 

\begin{theorem}\label{convergence}
  Let $\besov_R\subset\besov_D$ and
  let $u^\dagger$ be a minimum-$\norm{\cdot}_{\besov_R}$-solution of $Fu = v$.
  Then, for each minimizer $u^{\alpha,\delta}$ of
  \[
    T_\alpha(u) = \norm{Fu-v^\delta}^2 + \alpha \norm{u}_{\besov_R}^{p_R},
  \]
  and the parameter rule $\alpha \asymp \delta$ we get convergence 
  \begin{equation}\label{convergence_eq}
    \norm{u^{\alpha,\delta} - u^\dagger}_{\besov_R}
    \to 0, \quad \delta \to 0.
  \end{equation}
\end{theorem}

\begin{proof}[Proof of theorem~\ref{convergence}]
  We equip $\besov_D$ and $L_2$ with the weak topologies and want 
  to use theorem 3.5 from~\cite{hofmann2007convtikban}.
  To do so, we need the following to be fulfilled:
  \begin{enumerate}
  \item The norm $\norm{\cdot}_{L_2}$ is weakly lower-semicontinuous in $L_2$.
  \item $F:\besov_D\to L_2$ is weakly continuous.
  \item $\norm{\cdot}_{\besov_R}$ is proper, convex and weakly lower-semicontinuous on $\besov_D$.
  \item The sets $A_\alpha = \{u\ |\ \norm{Fu-v^\delta}^2 + \alpha \norm{u}_{\besov_R}^{p_R}<M\}$ are weakly sequentially compact in $\besov_D$.
  \end{enumerate}
  The first point is obvious and the second point is fulfilled by the assumption that $F$ is linear and continuous.
  For the forth point note that due to the continuous embedding we have
  $\norm{\cdot}_{\besov_D} \leq C \norm{\cdot}_{\besov_R}$ and hence,
  the sets $A_\alpha$ are bounded in $\besov_D$ which implies weak sequential compactness due to reflexivity of $\besov_D$.

  For the third point note that there exists a wavelet basis $\{\psi_\lambda\}$ which is an unconditional basis for both $\besov_R$ and $\besov_D$.
  Now, let $u_k\to u$ weakly in $\besov_D$ and $u_k\in\besov_R$.
  Since the $\psi_k$ are also elements of the dual spaces $\besov_R^*$ and $\besov_D^*$ it holds for all $\lambda$
  \[
  \ip{u_k}{\psi_\lambda}\stackrel{k\to\infty}{\longrightarrow}\ip{u}{\psi_\lambda}
  \]
  and hence, a sequence $u_k$ bounded in $\besov_R$ converges weakly to $u$ in $\besov_R$ if it does in the larger space $\besov_D$, because the duality pairing is the same in both spaces and $\{\psi_\lambda\}$ is an unconditional basis in $\besov_R$.
  This shows the weak lower-semicontinuity of $\norm{\cdot}_{\besov_R}$ on $\besov_R$-bounded sets in $\besov_D$ (which is sufficient for theorem 3.5 from~\cite{hofmann2007convtikban} to hold).

  Now, by theorem 3.5 from~\cite{hofmann2007convtikban} it follows, that $u^{\alpha,\delta}\to u^\dagger$ weakly in $\besov_D$ (and by the above considerations also in $\besov_R$) and moreover $\norm{u^{\alpha,\delta}}_{\besov_R} \to \norm{u^\dagger}_{\besov_R}$.
  Since $\besov_R$ is uniformly convex, this implies $u^{\alpha,\delta}\to u^\dagger$ strongly in $\besov_R$.
%
%
\end{proof}

Now we formulate a theorem on the rate of convergence which follows from the general results on regularization in Banach spaces~\cite{burger2004convarreg}.
We assume that certain knowledge on the true solution $u^\dagger$ is available,
i.e.~a certain source condition is fulfilled.
The source condition is formulated in terms of Besov smoothness.
This assumption, together with the assumptions on the range of $F^*$,
leads to a regularization term for which a certain convergence rate in Sobolev norm can be proven.

\begin{theorem}\label{convergence_rate}
  Let $u^\dagger \in \besov_S \subset \besov_D$ with $p_S \leq p_G$.
  Then, for each minimizer $u^{\alpha,\delta}$ of the Tikhonov functional
  \[
    T_\alpha(u) = \norm{Fu-v^\delta}^2 + \alpha \norm{u}_{\besov_R}^{p_R},
  \]
  with
  
  \begin{equation}
    \label{eq_def_pr_sr}
    p_R = \frac{p_S+p_G}{p_G},\quad
    s_R \leq \frac{p_S s_S - p_G s_G}{p_S + p_G}
  \end{equation}
  and the parameter rule $\alpha \asymp \delta$ we get the convergence rate
  \begin{equation}\label{convergence_rate_eq}
    \norm{u^{\alpha,\delta} - u^\dagger}_{H^{\sigma}}
    =
    \bigO(\sqrt{\delta}),
  \end{equation}
  where
  $\sigma: = s_R+d\Big(\frac{1}{2} - \frac{1}{p_R}\Big).$
\end{theorem}

\begin{remark}
  Notice, that in general the convergence statements in theorem \ref{convergence} and
  theorem \ref{convergence_rate} correspond to different Besov spaces
  $\besov_R$ and $H^{\sigma}$. The
  spaces coincide if and only if $p_S = p_G$, 
  otherwise we cannot give any information of inclusions,
  because the differential dimensions are equal:
  \[
    \sigma - \frac{d}{2} = s_R - \frac{d}{p_R}.
  \]
\end{remark}

\begin{remark}\label{remark_besovsource}
  The definitions of $p_R$ and $s_R$
  in theorem \ref{convergence_rate}
  imply that $\besov_S \subset \besov_R$.
  Otherwise the statement
  would not be meaningful, since if $\besov_S \supsetneqq \besov_R$,
  \[
    \exists \, u^\dagger \in \besov_S: \quad
    \norm{u^\dagger}_{\besov_R}^{p_R} = \infty.
  \]
  To see this note that due to $\besov_G \subset \besov_D^\ast$ the inequality 
  $\frac{1}{p_G} \geq  \frac{p_D -1}{p_D} $ holds and because
  $\besov_S \subset \besov_D$ we get 
  $\frac{1}{p_S} \geq \frac{1}{p_D}$
  and hence,
  \[
    p_R = \frac{p_S + p_G}{p_G} = p_S \Big(\frac{1}{p_G} + \frac{1}{p_S}\Big)
    \geq p_S \Big( \frac{p_D - 1}{p_D} + \frac{1}{p_D}\Big) = p_S.
  \]
  To see the inequality for the differential dimension of $\besov_R$ and $\besov_S$ note that $\besov_S \subset \besov_D \subset \besov_G^\ast$, and hence 
  \[
    -(s_S+s_G) < d\Big( \frac{p_G -1}{p_G} - \frac{1}{p_S} \Big),
  \]
  which leads to
  \begin{align*}
    -(s_S+s_G)p_Gp_S + & d(p_G+p_S) - d(p_Gp_S) \\
    < & d\Big( \frac{p_G -1}{p_G} - \frac{1}{p_S} \Big) p_Gp_S + d(p_G+p_S) - d(p_Gp_S)
        = 0.
  \end{align*}
  Applying this to the constraints (\ref{eq_def_pr_sr}) for $p_R$ and $s_R$ yields
  \begin{align*}
    s_R - \frac{d}{p_R} 
      & \leq \frac{p_S s_S - p_G s_G}{p_S + p_G} - d \frac{p_G}{p_S+p_G} \\
      & = s_S - \frac{d}{p_S} 
        + \frac{-(s_S+s_G)p_Gp_S + d(p_G+p_S) - d(p_Gp_S)}{(p_S+p_G)p_S}
        < s_S - \frac{d}{p_S}.
  \end {align*}
\end{remark}

For the proof of theorem \ref{convergence_rate} we need a property
of the mapping
\[
  \norm{\cdot}_{\besov_R}^{p_R}: \besov_S \to [0,\infty).
\]

\begin{proposition}\label{prop_besov_norm}
  Let $u \in \besov_S$ and let $s_R$ and $p_R$ fulfill~(\ref{eq_def_pr_sr}).
  Then
  \[
    \partial \Big(\norm{u}_{\besov_R}^{p_R} \Big) 
    = \Big\{ \nabla  \norm{u}_{\besov_R}^{p_R} \Big\}
    \subset \besov_G.
  \]
\end{proposition}

\begin{proof}
  Let $u \in \besov_S$.
  Since $p_S,p_G>0$, $p_R = 1 +  \frac{p_S}{p_G} > 1$, we get
  \begin{equation}\label{subgradient_eq}
    \partial(\norm{u}_{\besov_R}^{p_R}) 
    = \{ \nabla \norm{u}_{\besov_R}^{p_R} \}
    = \Big\{ p_R \; \sum\limits_{\lambda\in\Lambda} 
      2^{p_R(s_R+d(\frac{1}{2}-\frac{1}{p_R})) \; |\lambda|}
      \sign (u_\lambda) \; |u_\lambda|^{p_R-1}\Big\},
  \end{equation}
  and hence,
  \begin{align*}
    \Big\| \nabla \norm{u}_{\besov_R}^{p_R} \Big\|_{\besov_G}^{p_G}
    & = \sum\limits_{\lambda\in\Lambda} 
      2^{p_G(s_G+d(\frac{1}{2}-\frac{1}{p_G})) \; |\lambda|}
      \Big| 
      p_R \;
      2^{p_R(s_R+d(\frac{1}{2}-\frac{1}{p_R})) \; |\lambda|}
      \sign (u_\lambda) \; |u_\lambda|^{p_R-1}
      \Big|^{p_G}\\
    & = p_R^{p_G} \; \sum\limits_{\lambda\in\Lambda} 
      2^{[p_G(s_G+d(\frac{1}{2}-\frac{1}{p_G}))
        +p_G p_R(s_R+d(\frac{1}{2}-\frac{1}{p_R}))]
        \; |\lambda|}
      |u_\lambda|^{p_G(p_R-1)}.
  \end{align*}
  Because~(\ref{eq_def_pr_sr}) if fulfilled we get
  for the exponent
  \begin{align*}
      & p_G\Big(s_G+d\Big(\frac{1}{2}-\frac{1}{p_G}\Big)\Big)
        +p_G p_R\Big(s_R+d\Big(\frac{1}{2}-\frac{1}{p_R}\Big)\Big)\\
    = & p_G s_G + p_G p_R s_R + d\Big(\frac{p_G p_R}{2} - \frac{p_G}{2} - 1\Big) \\
    = & p_S \Big(\frac{p_G}{p_S} s_G + \frac{p_G}{p_S} s_R + s_R\Big)
        + d p_S \Big( \frac{1}{2} - \frac{1}{p_S} \Big)
        \leq p_S s_S + p_S d \Big(\frac{1}{2} - \frac{1}{p_S} \Big),
  \end{align*}
  hence
  \[
    \Big\| \nabla \norm{u}_{\besov_R}^{p_R} \Big\|_{\besov_G}^{p_G}
    \leq p_R^{p_G} \; \sum 2^{p_S(s_S + d(\frac{1}{2} - \frac{1}{p_S})) \; |\lambda|} |u_\lambda|^{p_S} 
    \asymp p_R^{p_G} \norm{u}_{\besov_S}^{p_S}
    < \infty,
  \]
  since $u\in\besov_S$.
\end{proof}

Now we are able to do the
\begin{proof}[Proof of theorem \ref{convergence_rate}]
  In \cite{burger2004convarreg} it is proved that the source condition
  \begin{equation}\label{source_cond}
    \exists \, w \in L_2: \quad
    F^\ast w \in \partial(\norm{u^\dagger}_{\besov_R}^{p_R})
  \end{equation}
  leads to the estimate for the so-called Bregman distance
  \[
    D_{\partial \norm{u^\dagger}_{\besov_R}^{p_R}} (u^{\alpha,\delta},u^\dagger) = \bigO(\delta)
  \]
  for minimizers $u^{\alpha,\delta}$ of the Tikhonov functional~(\ref{tikhonov_eq}) and $\alpha \asymp \delta$.
  Here by assumption the range of the adjoint operator $F^\ast$ is
  $\besov_G$,
  and hence, with proposition~\ref{prop_besov_norm}, we get
  \[
    \partial \Big(\norm{u}_{\besov_R}^{p_R} \Big) 
    = \Big\{ \nabla  \norm{u}_{\besov_R}^{p_R} \Big\}
    \subset \besov_G
    = \rg (F^\ast),
  \]
  thus the source condition (\ref{source_cond}) is fulfilled.
  Further we get with (\ref{subgradient_eq})
  \begin{align*}
    D_{\nabla \norm{u^\dagger}_{\besov_R}^{p_R}} & (u^{\alpha,\delta},u^\dagger)
    = \norm{u^{\alpha,\delta}}_{\besov_R}^{p_R}
      -\norm{u^\dagger}_{\besov_R}^{p_R}
      - \ip{\nabla \norm{u^\dagger}_{\besov_R}^{p_R}}{u^{\alpha,\delta}-u^\dagger} \\
    = & \sum\limits_{\lambda\in\Lambda}
      2^{p_R(s_R+d(\frac{1}{2} - \frac{1}{p_R})) \; |\lambda|}
      \Big(
      |u^{\alpha,\delta}_\lambda|^{p_R}
      -|u^\dagger_\lambda|^{p_R} \\
      & - p_R \sign (u^\dagger_\lambda) \; |u^\dagger_\lambda|^{p_R-1}
        (u^{\alpha,\delta}_\lambda - u^\dagger_\lambda)\Big)\\
    = & \bigO(\delta).
  \end{align*}
  For $a,b \in \R$, $C>|a|$, $|b-a|<L$, $1<p\leq 2$ by \cite[lemma 4.7]{bredies2008harditer},
  \[
    |b|^p - |a|^p - p \sign (a) |a|^{p-1} (b-a) \geq k(p,C,L) |b-a|^2,
  \]
  where $k(p,C,L)$ is a positive constant which depends on $p$, $C$ and $L$.
  Since by remark~\ref{remark_besovsource} it holds $u^\dagger \in \besov_S \subset \besov_R$ and hence,
  \[
    \exists \, C > 0: \quad
    \Big|
    2^{(s_R+d(\frac{1}{2}-\frac{1}{p_R})) \; |\lambda|}
    |u^\dagger_\lambda|
    \Big|<C, \quad \lambda \in \Lambda.
  \]
  Furthermore, since $\norm{u^{\alpha,\delta} - u^\dagger}_{\besov_R} \to 0$ for $\delta \to 0$ 
  we get according to theorem \ref{convergence}
  \[
    \exists \, L > 0: \quad
    \Big|
    2^{(s_R+d(\frac{1}{2}-\frac{1}{p_R})) \; |\lambda|}
    [u^{\alpha,\delta}_\lambda - u^\dagger_\lambda]
    \Big|<L, \quad \lambda \in \Lambda.
  \]
  Applying this with 
  $a=2^{(s_R+d(\frac{1}{2}-\frac{1}{p_R})) \; |\lambda|} u^{\alpha,\delta}_\lambda$, 
  $b=2^{(s_R+d(\frac{1}{2}-\frac{1}{p_R})) \; |\lambda|} u^\dagger_\lambda$,
  and $p=p_R \in (1,2]$, since $p_G \geq p_S$,
  we get
  \begin{align*}
    D_{\nabla \norm{u^\dagger}_{\besov_R}^{p_R}} & (u^{\alpha,\delta},u^\dagger)
      \geq K \sum\limits_{\lambda\in\Lambda} 
      2^{(2(s_R+d(\frac{1}{2}-\frac{1}{p_R}))) \; |\lambda|} 
      |u^{\alpha,\delta}_\lambda - u^\dagger_\lambda|^2\\
    & \asymp K 
      \norm{u^{\alpha,\delta}_\lambda - u^\dagger_\lambda}_{H^{s_R+d\big(\frac{1}{2} - \frac{1}{p_R}\big)}}^2,
  \end{align*}
  because of the norm equivalence~(\ref{eq_equiv_besov_norm}) and the fact that $H^s=B^s_{2,2}$ for all $s$.
  
  Finally, this gives
  \[
  \norm{u^{\alpha,\delta} - u^\dagger}_{H^{\sigma}}
  =
  \bigO(\sqrt{\delta})
  \]
  where $\sigma: = s_R+d\Big(\frac{1}{2} - \frac{1}{p_R}\Big)$.
\end{proof}


\section{Source Condition Weakening}\label{sec_optimization}
In the setup of theorem~\ref{convergence_rate} we assumed that a source condition in terms of 
Besov smoothness is known, i.e.~$u^\dagger\in\besov_S$.
From that a regularization penalty $\norm{\cdot}_{\besov_R}^{p_R}$ was derived
which leads to a convergence rate in a certain Sobolev space.

\begin{figure}[h]
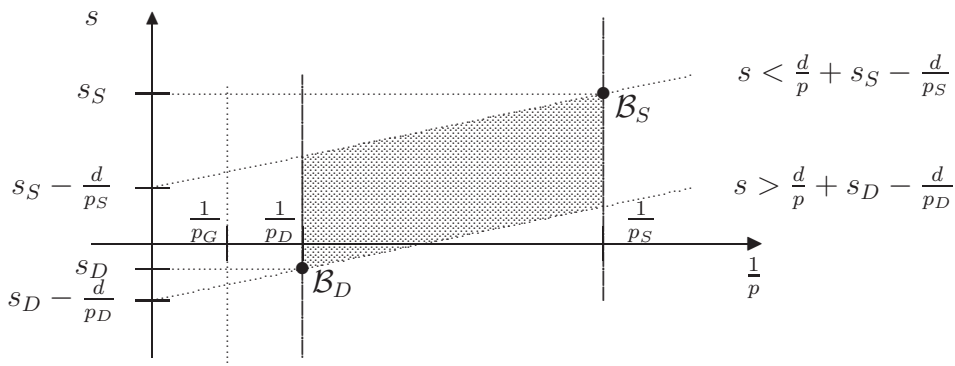

  \[
  \beginpicture
    \setcoordinatesystem units <4cm,1.5cm>
    \setplotarea x from -0.2 to 2, y from -1 to 2
    \axis left shiftedto x=0
          ticks long in unlabeled at -0.5 0.5 -0.222222 1.333333 /
                long out unlabeled at -0.5 0.5 -0.222222 1.333333 /
    /
    \axis bottom shiftedto y=0
          ticks long in unlabeled at 0.25 0.5  1.5 /
                long out unlabeled at 0.25 0.5 1.5 /  
    /
    \put {$\scriptstyle \blacktriangle$} at 0 2
    \put{$s$} at -0.2 2
    \put{$\scriptstyle \blacktriangleright$} at 2 0
    \put{$\frac{1}{p}$} at 2 -0.3
    
    \setdots <2pt>
    \plot 0 -0.5 1.8 0.5 /
    \put{$s > \frac{d}{p} + s_D - \frac{d}{p_D}$} at 2.3 0.5
    \put{$s_D - \frac{d}{p_D}$} at -0.3 -0.5
    \plot 0  0.5 1.8 1.5 /
    \put{$s < \frac{d}{p} + s_S - \frac{d}{p_S}$} at 2.3 1.5
    \put{$s_S - \frac{d}{p_S}$} at -0.3 0.5
    
    \plot 0 -0.222222  0.5 -0.222222 /
    \put{$s_D$} at -0.2 -0.222222
    \plot 0 1.333333  1.5 1.333333 /
    \put{$s_S$} at -0.2 1.333333
    
    \put{$\bullet$} at 1.5 1.333333
    \put{$\besov_S$} at 1.6 1.2
    
    \put{$\bullet$} at 0.5 -0.222222
    \put{$\besov_D$} at 0.6 -0.35
        
    \plot 0.25 1.3889  0.25 -1.1111 /
    
    \setsolid
    \plot 0.5 1.5  0.5 -1 /
    \plot 1.5 2  1.5 -0.5 /
    
    \put{$\frac{1}{p_G}$} at 0.175 0.2
    \put{$\frac{1}{p_D}$} at 0.425 0.2    
    \put{$\frac{1}{p_S}$} at 1.625 0.2
    
    \setshadegrid span <1pt>
    \vshade 0.5 -0.222222 0.777777
            <z,z,z,z>
            1.5  0.333333 1.333333
    /
  \endpicture
  \]
  \caption{Weaker source conditions $u^\dagger \in B_p^s$ 
    with $\besov_S \subset B_p^s \subset \besov_D$, $\frac{1}{p_D} > \frac{1}{p_G}$}
  \label{bild_feas_reg}
\end{figure}

Besov spaces are embedded into each other 
via the non-linear intricate properties (\ref{eq-embed_besov_spaces}) of proposition~\ref{embed_besov_spaces}.
Considering this, the question arises which penalties $\norm{\cdot}_{\besov_R}^{p_R}$
and convergence rates (i.e.~which $\sigma$) follows from a weakened source condition
$u^\dagger \in B_p^s$
with $\besov_S \subset B_p^s \subset \besov_D$.
In addition to the embedding properties (\ref{eq-embed_besov_spaces}) for application of
theorem~\ref{convergence_rate} one has to ensure $p \leq p_G$. This yields to
the following set of possible weaker source conditions
  \begin{align}
    u^\dagger \in B_p^s\quad\text{such that}\quad
      & \besov_S \subset B_p^s \text{, i.e.} 
        & s_S - \frac{d}{p_S} & > s - \frac{d}{p}, \label{feas_reg1}\\
      & & \frac{1}{p_S} & \geq \frac{1}{p}, \label{feas_reg2}\\
      & B_p^s \subset \besov_D \text{, i.e.} 
        & s - \frac{d}{p} & > s_D - \frac{d}{p_D}, \label{feas_reg3}\\
      & & \frac{1}{p} & \geq \frac{1}{p_D}, \label{feas_reg4}\\
      & p_R \in (1,2] \text{, i.e.}
        & \frac{1}{p} & \geq \frac{1}{p_G}. \label{feas_reg5}
  \end{align}
Figure~\ref{bild_feas_reg} illustrates the set of weaker source conditions,
i.e. the equalities~(\ref{feas_reg1})-(\ref{feas_reg5}), graphically
for $p_D<p_G$, which ensures $p\leq p_G$.

The direct application of theorem~\ref{convergence_rate} to the idea of weakening the source condition
with (\ref{feas_reg1})-(\ref{feas_reg5}) gives the following theorem.
\begin{theorem}\label{convergence_rate_weaker_source}
  Let $u^\dagger \in \besov_S \subset \besov_D$ with $p_S \leq p_G$.
  Further let $p>0$ with $p_S \leq p \leq \min \{p_D,p_G\}$ and $\varepsilon>0$.
  Then, for each minimizer $u^{\alpha,\delta}$ of the Tikhonov functional
  \[
    T_\alpha(u) = \norm{Fu-v^\delta}^2 + \alpha \norm{u}_{\besov_R}^{p_R},
  \]
  with
  \begin{equation}\label{eq_def_pr_sr_weaker_source}
    p_R = \frac{p+p_G}{p_G},\quad
    s_R \leq \frac{p \, s_S - p_G s_G}{p + p_G} - d \Big(\frac{1}{p + p_G} \Big(\frac{p}{p_S} - 1\Big)\Big) 
      - \varepsilon \Big(\frac{p}{p_S} - 1\Big),
  \end{equation}
  and the parameter rule $\alpha \asymp \delta$ we get the convergence rate
  \[
    \norm{u^{\alpha,\delta} - u^\dagger}_{H^{\sigma}}
    =
    \bigO(\sqrt{\delta}),
  \]
  where
  $\sigma: = s_R+d\Big(\frac{1}{2} - \frac{1}{p_R}\Big).$
\end{theorem}

\begin{proof}
  From (\ref{feas_reg2}),(\ref{feas_reg4}) and (\ref{feas_reg5}) it follows that theorem~\ref{convergence_rate}
  is applicable for $p>0$ with $p_S \leq p \leq \min \{p_D,p_G\}$. To ensure the embedding properties 
  for the differential dimension, i.e. equation (\ref{feas_reg1}) and (\ref{feas_reg3}), 
  one has to choose $s\in\R$ with
  \[
    s_S - \frac{d}{p_S} > s - \frac{d}{p} > s_D - \frac{d}{p_D}.
  \]
  With that the application of theorem~\ref{convergence_rate} yields (\ref{eq_def_pr_sr_weaker_source})
  and hence the convergence in $H^\sigma$.
\end{proof}

The convergence result in theorem~\ref{convergence_rate_weaker_source} gets stronger
as $\sigma$ increases.
Since $\sigma$ depends on $s_R$ and $p_R$ we address the question
how to choose $\besov_R$ in a way such that $\sigma$ is maximal.
We try to find the regularization penalty $\norm{\cdot}_{\besov_R}^{p_R}$,
which gives the best estimate with respect to $\sigma$
(while $F$ and the spaces $\besov_S$, $\besov_G$ and $\besov_D$ are fixed).
Since $\sigma$ depends strictly monotone on $s_R$, we have to choose 
$s_R$ as large as possible so that we have to solve the following
optimization problem, which depends only on $p$:
\begin{equation}
  \label{eq_opt_problems}
  \left.
    \begin{array}{rl}
      \max\limits_{p} \quad &  
      \frac{p s_S - p_G s_G}{p + p_G} + d 
      \Big[\frac{1}{2} - \frac{1}{p + p_G} 
      \Big(\frac{p_G}{p_G^\ast} + \frac{p}{p_S} \Big)
      \Big] - \varepsilon \Big(\frac{p}{p_S} - 1\Big)\\
      \text{such that}\quad
      & p_S \leq p \leq \min \{p_D,p_G\}.
    \end{array}
    \right\}
\end{equation}

Since $\varepsilon>0$ can be arbitrarily small we neglect the term 
$\varepsilon \big(\frac{p}{p_S} - 1\big)$ and hence, we have to find the 
maximum of
\begin{align*}
  \widehat{\sigma}(p)
  & := \frac{p s_S - p_G s_G}{p + p_G} + d \Big[\frac{1}{2} - \frac{1}{p + p_G} 
       \Big(\frac{p_G}{p_G^\ast} + \frac{p}{p_S} \Big)\Big]\\
  & = \frac{p}{p+p_G} \Big(s_S - \frac{d}{p_S} \Big)
    + \frac{p_G}{p+p_G} \Big(-s_G - \frac{d}{p_G^\ast} \Big).
\end{align*}
The function $\widehat{\sigma}$ is monotonically increasing in $p$, since for $p_1>p_2$ we get
\begin{align*}
  \widehat{\sigma}&(p_1) -\widehat{\sigma}(p_2)\\
  & = \Big(\frac{p_1}{p_1+p_G}-\frac{p_1}{p_1+p_G}\Big) \Big(s_S - \frac{d}{p_S} \Big)
    + \Big(\frac{p_G}{p_2+p_G}-\frac{p_G}{p_2+p_G}\Big) \Big(-s_G - \frac{d}{p_G^\ast}\Big)\\
  & = \Big(p_G \frac{p_1-p_2}{(p_1+p_G)(p_2+p_G)} \Big)
      \Big(s_S - \frac{d}{p_S} - \Big(-s_G - \frac{d}{p_G^\ast}\Big) \Big) > 0,
\end{align*}
since $\besov_S \subset \besov_D \subset \besov_G^\ast$.
This proves
\begin{corollar}\label{opt_convergence_rate}
  Let $u^\dagger \in \besov_S \subset \besov_D$, $p_S \leq p_G$
  and $\varepsilon>0$ sufficiently small.
  Then the Tikhonov regularization $T_\alpha$
  with the parameter rule $\alpha \asymp \delta$,
  penalty according to (\ref{eq_def_pr_sr_weaker_source}) 
  with 
  \[
    p:=min\{p_D,p_G\}
  \]
  gives the strongest convergence.
  \begin{enumerate}[i)]
    \item If $p_G \geq p_D$, we get with $\widetilde{\varepsilon}:=\varepsilon \Big(\frac{p}{p_S} - 1\Big)$
      \begin{align*}
        p_R & = \frac{p_D+p_G}{p_G},\\
        s_R & = \frac{p_D s_S - p_G s_G}{p_D + p_G} - d (\frac{1}{p_D + p_G} (\frac{p_D}{p_S} - 1)) - \widetilde{\varepsilon},\\
      \intertext{a convergence rate result (\ref{convergence_rate_eq}) in $H^\sigma$ with}
        \sigma & = \frac{p_D s_S - p_G s_G}{p_D + p_G} + d 
        \Big[\frac{1}{2} - \frac{1}{p_D + p_G} 
        \Big(\frac{p_G}{p_G^\ast} + \frac{p_D}{p_S} \Big)
        \Big] - \widetilde{\varepsilon}.
      \end{align*}
    \item If $p_G < p_D$, we obtain with
      \begin{align*}
        p_R & = 2,\\
        s_R & = \frac{1}{2} \big(s_S - \frac{d}{p_S} - \big(s_G - \frac{d}{p_G}\big)\big)- \widetilde{\varepsilon},
      \end{align*}
      convergence rate~(\ref{convergence_rate_eq}) in $\besov_R = H^{s_R}$.
  \end{enumerate}
\end{corollar}

\begin{remark}
  \label{rem:no_optimization_possible}
  If $p_S < \min \{p_D, p_G\}$ the convergence rate in corollary~\ref{opt_convergence_rate} is better than in 
  theorem~\ref{convergence_rate}. Note that the Tikhonov functionals do not
  coincide.
\end{remark}

A curiosity of theorem~\ref{convergence_rate}, i.e. 
of the straight forward application of the Banach space regularization 
results~\cite{burger2004convarreg,hofmann2007convtikban},
is that a more restrictive
source condition $\besov_T \subset \besov_S$ 
does not necessarily enforce a better convergence rate.
As the following counterexample shows, sometimes the converse may happen.
\begin{counterexample}\label{counterex}
  Let $\varepsilon > 0$ be sufficiently small and $\eta >0$. Further let $F$ be an operator with
  \[
    F: H^{-\eta} \to L_2, \quad
    \rg F^\ast = H^{\eta}.
  \]
  \begin{enumerate}
    \item With the loose source condition $u^\dagger \in \besov_S = H^\eta$
      theorem~\ref{convergence_rate} yields a convergence rate in Lebesgue space
      $L_2$.
      (The choice $\besov_S = H^\eta$ leads to $\besov_R = L^2$.)
    \item If we tighten the condition to 
      $u^\dagger \in B_1^{\eta + \frac{d}{2} + 3\varepsilon} \subset H^\eta$
      we get the regularization space $\besov_R = B^{\eta+\varepsilon+d/6}_{3/2}$
      and a convergence rate in Sobolev space $H_2^{-\frac{\eta}{3}+\varepsilon}$,
      which is larger than in $L_2$ for small $\varepsilon$.
  \end{enumerate}
\end{counterexample}

In contrast to that the usage of Besov space embeddings, i.e.~corollary~\ref{opt_convergence_rate}, rewards a 
tighter source condition with a stronger convergence rate:
Let $\besov_T \subset \besov_S$, i.e.
\[
  s_T - \frac{d}{p_T} > s_S - \frac{d}{p_S},
  \quad
  p_T \leq p_S.
\]
Then we get consequently for case \textit{i)} ($p_G\geq p_D$)
\begin{equation*}
  \sigma(\besov_T) - \sigma(\besov_S) 
  = \frac{p_D}{p_D+p_G} \Big( s_T - \frac{d}{p_T} - \big( s_S - \frac{d}{p_S}\big)\Big) > 0,
\end{equation*}
and for case \textit{ii)} ($p_G < p_D$)
\begin{equation*}
  \sigma(\besov_T) - \sigma(\besov_S) 
   = \frac{1}{2} \Big( s_T - \frac{d}{p_T} - \big( s_S - \frac{d}{p_S}\big)\Big) > 0.
\end{equation*}


\section{Examples}
\label{sec_examples}
In the following we will illustrate the convergence rate results with a few examples.
With the first one we want to show that the choice of the parameter $p$ resp. 
the choice of the source condition $\besov_S \subset B_p^s \subset \besov_D$
(cf.~(\ref{feas_reg1})-(\ref{feas_reg5})) 
influences the convergence rate significantly.

\begin{example}[Smoothing in the Sobolev scale]
  Let $d=1$, $\eta>\frac{1}{2}$ and consider the operator
  \[
    F: H^{-\eta} \to L_2, \quad \text{ with } \rg F^\ast = H^\eta,
  \]
  i.e.~we consider smoothing of order $\eta$ in the Sobolev scale.
  Moreover, we assume that the source condition $u^\dagger \in B_1^{2\eta} \subset H^{-\eta}$ holds.
  Theorem~\ref{convergence_rate_weaker_source} yields convergence rates for Tikhonov penalties
  $\norm{\cdot}_{B_{p_R(p)}^{s_R(p)}}^{p_R(p)}$ with $p_S \leq p \leq p_G = p_D$.
  Since $\sigma$ s monotone in $p$, cf.~solution of the optimization problem, we just investigate in the two boundary values here.
  For $p=p_S$ we get the Tikhonov functional
  \[
    T_\alpha(u) = \norm{Fu - v^\delta}^2 + \alpha \norm{u}_{B_{\frac{3}{2}}^0}^{\frac{3}{2}}.
  \]
  With that worst parameter choice resp. worst source condition
  theorem \ref{convergence_rate_weaker_source} yields 
  \[
    \sigma = \frac{2 \eta - 2 \eta}{3} + \frac{1}{2} - \frac{2}{3} = - \frac{1}{6}
  \]
  and hence, the convergence rate occurs in a Sobolev space $H^\sigma$ with negative smoothness.
  
  Next let us check the rate with optimal parameter $p=p_G$. Here we get
  for the Tikhonov functional
  \[
    T_\alpha(u) = \norm{Fu - v^\delta}^2 + \alpha \norm{u}_{H^{\frac{\eta}{2} - \frac{1}{4} - \varepsilon}}^2
  \]
  the convergence rate in the Sobolev space $H^\sigma$ with smoothness
  \[
    \sigma = \frac{4 \eta - 2 \eta}{4} + \frac{1}{2} - \frac{1}{4} (1 + 2) -\varepsilon
           = \frac{\eta}{2} - \frac{1}{4} -\varepsilon,
  \]
  which is greater than zero for small $\varepsilon$, since $\eta > \frac{1}{2}$. 
  Hence, we get a convergence rate in a Sobolev space with positive smoothness. 
\end{example}

The first example may lead to the conclusion that a penalty
formulated in a Sobolev space gives the best convergence rate.
This  impression may be intensified, because the optimal source 
also lives in a Sobolev spaces, i.e. $p=p_D=2$.
As we will see now with the next two examples with operators formulated in Banach scales,
this guess is not true.
Moreover, the following examples illustrate the difference
between the cases $p_S = \min\{p_D,p_G\}$ and $p_S < \min\{p_D,p_G\}$.
In the first case theorem~\ref{convergence_rate_weaker_source} yields a convergence rate
for only one Tikhonov functional resp. no optimization is possible, cf. example~\ref{ex_Banach1}.
In the second case we get a set of allowed Tikhonov penalties depending on $p$,
$p_S \leq p \leq \min\{p_D,p_G\}$.

\begin{example}[Smoothing in the Besov scale, $p_S = \min\{p_D,p_G\}$]\label{ex_Banach1}
  Let $d=1$, $\eta > 0$ arbitrary, $0 < \theta \leq 1$ small and consider the operator
  \[
    F:B_{1+\theta}^{-\eta} \to L_2, \quad \text{ with } \rg F^\ast = B_{\frac{1+\theta}{\theta}}^\eta = \Big(B_{1+\theta}^{-\eta}\Big)^\ast,
  \]
  which models smoothing in the scale of Besov spaces.
  Further let $u^\dagger \in B_{1+\theta}^{-\eta+\theta}$ be the source condition.
  Notice that $B_{1+\theta}^{-\eta+\theta} \subset B_{1+\theta}^{-\eta}$
  and 
  \[
    1+\theta \leq \frac{1+\theta}{\theta},
    \quad\mbox{for } \theta \leq 1,
  \]
  and hence, we can guarantee $p_S \leq p_G$.
  Due to $p_S = \min\{p_D,p_G\} = p_D$ it follows from theorem~\ref{convergence_rate_weaker_source}
  that only the Tikhonov functionals with $p=p_S=p_D$, i.e. 
  \[
    T_\alpha(u) = \norm{Fu - v^\delta}^2 + \alpha \norm{u}_{\besov_R}^{p_R},
  \]
  with $p_R = \frac{p+p^\ast}{p^\ast} = p = \theta+1$ and
  \[
    s_R \leq \frac{p (-\eta+\theta) + p^\ast (-\eta)}{p + p^\ast}
        = -\eta + \frac{\theta^2}{\theta + 1}
  \]
  yields a convergence rate in Sobolev space $H^\sigma$.
  The maximal smoothness $\sigma$ is obtained with the penalty with $s_R = -\eta + \frac{\theta}{p^\ast}$ and it reads as
  \[
    \sigma 
    = \frac{p (-\eta+\theta) + p^\ast (-\eta)}{p + p^\ast} 
    + \frac{1}{2}-\frac{1}{p}
    = -\eta + \frac{1}{2} - \frac{\theta - 1}{\theta + 1}.
  \]

  \emph{To put it roughly}: For the operator
  \[
    F: B_1^{-\eta} \to L_2 \quad \text{ with } \rg F^\ast = B_\infty^{\eta},
  \]
  the source condition $u^\dagger \in B_1^{-\eta}$ and
  the Tikhonov functional
  \[
    T_\alpha(u) = \norm{Fu - v^\delta}^2 + \alpha \norm{u}_{B_1^{-\eta}}^{1},
  \]
  we a get convergence rate
  \[
    \norm{u^{\alpha,\delta} - u^\dagger}_{H^{-\eta - \frac{1}{2}}}
    = \bigO(\sqrt{\delta}).
  \]
\end{example}

In the above Besov scale example no optimization of the convergence rate was possible. In the next example
there is a set of allowed Tikhonov regularizations and hence an optimal one.

\begin{example}[Smoothing in the Besov scale, $p_S < \min\{p_D,p_G\}$]\label{ex_Banach2}
   Let $d=1$, $\eta > 0$ arbitrary, $0 < \theta \leq \frac{1}{2}$ small and consider the operator
  \[
    F:B_{\frac{3}{2}}^{-\eta} \to L_2, \quad \text{ with } \rg F^\ast = B_{3}^\eta = \Big(B_{\frac{3}{2}}^{-\eta}\Big)^\ast.
  \]
  Further let $u^\dagger \in B_{1+\theta}^{-\eta+1} \subset B_{\frac{3}{2}}^{-\eta}$ be the source condition.
  Notice that since $\besov_S \subset \besov_D$ and $p_D < p_G$
  we can guarantee the second assumption of theorem~\ref{convergence_rate_weaker_source},
  $p_S \leq p_G$.
  Here theorem~\ref{convergence_rate_weaker_source} yields convergence rates for a set of Tikhonov penalties
  with $p_S \leq p \leq p_D$.
  We will just investigate in the two boundary values here again.

  With the worst parameter choice, i.e. $p=p_S$,
  theorem \ref{convergence_rate_weaker_source} yields with the Tikhonov functional
  \[
    T_\alpha(u) = \norm{Fu - v^\delta}^2 + \alpha \norm{u}_{\besov_R}^{p_R},
  \]
  with $p_R = \frac{4+\theta}{3}$ and
  $s_R = \frac{p_S s_S - p_G s_G}{p_S + p_G} = -\eta + \frac{\theta + 1}{\theta +4}$
  a convergence rate in $H^\sigma$ with
  \[
    \sigma = -\eta + \frac{1}{2} + \frac{\theta - 2}{\theta + 4}.
  \]

  For the optimal parameter $p=p_D$ we get a penalty $p_R = p_D = \frac{3}{2}$ and 
  \[
    s_R = \frac{1}{p_G} s_S - \frac{1}{p_D} s_G - \Big(\frac{1}{p_D+p_G} \Big(\frac{p_D}{p_S} - 1\Big) \Big) - \widetilde\varepsilon
    = -\eta + \frac{1}{3} + \frac{1}{9} \Big(\frac{2\theta - 1}{\theta + 1}\Big) - \widetilde\varepsilon.
  \]
  Theorem~\ref{convergence_rate_weaker_source} yields a convergence rate with
  \[
    \sigma 
    = -\eta + \frac{1}{6} + \frac{1}{9} \Big(\frac{2\theta - 1}{\theta + 1}\Big) - \widetilde\varepsilon
  \]
  with small $\widetilde\varepsilon:=\varepsilon \Big(\frac{p_D}{p_S} - 1\Big)>0$.

  \emph{To put it roughly}: Consider the operator
  \[
    F:B_{\frac{3}{2}}^{-\eta} \to L_2,\quad
    F^\ast: L_2 \to B_{3}^\eta,
  \]
  and the source condition $u^\dagger \in B_1^{-\eta+1}$.

  For the worst choice $p=p_S$ we get for the Tikhonov functional
  \[
    T_\alpha(u) = \norm{Fu - v^\delta}^2 + \alpha \norm{u}_{B_{\frac{4}{3}}^{-\eta+\frac{1}{4}}}^{\frac{4}{3}},
  \]
  a convergence rate
  \[
    \norm{u^{\alpha,\delta} - u^\dagger}_{H^{-\eta}}
    = \bigO(\sqrt{\delta}).
  \]

  The optimal choice $p=p_D$ yields with the Tikhonov functional
  \[
    T_\alpha(u) = \norm{Fu - v^\delta}^2 + \alpha \norm{u}_{B_{\frac{3}{2}}^{-\eta+\frac{2}{9}}}^{\frac{3}{2}},
  \]
  a convergence rate
  \[
    \norm{u^{\alpha,\delta} - u^\dagger}_{H^{-\eta+\frac{1}{18}}}
    = \bigO(\sqrt{\delta}).
  \]
\end{example}


\section{Conclusion}
The aim of this paper was to make a first attempt to analyze scales of Banach spaces for Tikhonov regularization.
We used Besov spaces to model the smoothing properties of the operator, the regularization term and the source condition.
The convergence rates results were obtained in the Hilbert scale of Sobolev spaces.
In comparison to regularization in Hilbert scales initiated in~\cite{natterer1984hilbertscales} the relation between these spaces is more complicated.
Of particular interest is the fact that on the one hand tighter source conditions may not lead to stronger convergence rates and on the other hand a less tight source condition may result in a stronger result.

Our examples in section~\ref{sec_examples} show only slight improvements in the Sobolev exponents when the Besov-penalty is optimized.
It is questionable if the effect can be observed numerically.
However, the effect that looser source conditions lead to tighter convergence results is interesting on its own

We did not use Besov spaces neither for the discrepancy term in the Tikhonov functional nor to measure the convergence rate.
Both points are of interest and may lead to more general results.
Since this paper is a first attempt in this direction we postpone this analysis for future work.


\section*{References}
\bibliographystyle{plain}
\bibliography{literature}

\end{document}